\documentclass[12pt]{article}

\usepackage{graphicx}
\usepackage{amssymb}
\usepackage{amsfonts}

\usepackage{times,amsmath,amsthm,amssymb,latexsym,amscd}

\newtheorem{Th}{Theorem}[section]

\newcommand{\btbt}{\left( \begin{array}{cc}}
\newcommand{\etbt}{ \end{array}\right)}
\newcommand{\bthbth}{\left( \begin{array}{ccc}}
\newcommand{\ethbth}{ \end{array}\right)}
\newtheorem{Lm}{Lemma}[section]

\newtheorem{co}{Corollary}[section]

\newcommand{\bcol}{\left(\begin{array}{c}}
\newcommand{\ecol}{\end{array}\right)}

\title{Regular trees in random regular graphs}
\author{Eran Makover\thanks{makovere@ccsu.edu} and Jeffrey McGowan\thanks{mcgowan@ccsu.edu}}
\begin{document}

\maketitle
\begin{abstract}We investigate the size of the embedded regular tree rooted at a vertex in a $d$ regular random graph.  We show that almost always, the size of this tree will be $\frac{1}{2}\log n$, where $n$ is the number of vertices in the graph.  We use this to give an asymptotic estimate for Gauss' Hypergeometric Function.\end{abstract}
\section{Introduction} 
Random regular graphs turn out to be quite different from the usual models of random graphs, and fewer results exist on their basic properties \cite{bollobas}  \cite{janson}. In this short note we use the configuration model for random regular graphs with $n$ vertices.


We  pick a vertex with uniform distribution and we estimate the size of the embedded  regular tree rooted in the given vertex.  
Our main result is that for $d$ fixed, in a $d$ regular graph with $n$ vertices, the expected value for the size of the $d$ regular embedded tree rooted at a given vertex is almost always $\sqrt{n}$

\begin{Th}\label{main}
Let $G$ be a random $d$-regular graph with $|V|=n$. For $v\in G$ the radius of the largest embedded $d$-regular tree rooted in $v$ will almost always be $(\frac{1}{2}+o(1))\log_{d-1} n$
\end{Th}

Define the radius of a rooted tree as $r=\max \delta(root,v)$ where $\delta$ is the distance in the graph.  Then for $d$-regular rooted tree with $n$ vertices,  $r\propto \log_{d-1} n$, our result shows that $r= (\frac{1}{2} +o(1))\ \log_{d-1} n$ almost always for the embedded tree in a $d$ regular graph. We note that since almost all  random regular graphs  are expanding graphs,  their diameter  is $O(\log n)$ \cite{ HLW}, and therefore for every vertex we can find a spanning tree rooted in the vertex with $r=O(\log n)$.

Our initial motivation for considering this question was a desire  to find fundamental domains for Belyi surfaces.  Such surfaces come from oriented three regular graphs \cite{Belyi, Groth}, and spanning trees in these graphs correspond to fundamental domains on the open surfaces \cite{Rivin}. 
We treat the problem for a graph in a similar way to finding a Dirichlet fundamental domain for a surface.  Pick a random vertex (root) in the graph and construct a spanning tree by connecting each vertex too the root by the shortest path.  If more then one path exists between the root and a given vertex choose either of the paths and use the same path to all points which are farther from the root.   We move "layer by layer" out from the root, and consider when some vertex in the outermost layer is connected to an already traversed vertex.  Intuitively, this is essentially the birthday problem; how big a party will make it likely that there will be two people who have the same birthday?

Theorem \ref{main} shows that for $\frac{1}{2}$ of the radius this "Dirichlet spanning tree" is regular. It is well known the expander graphs are locally like trees, but here we give a more precise description of this situation.  Our techniques are related to those used by Bollobas and de la Vega \cite{BV}, They investigated the diameter of random regular graphs, but did not consider the size of a rooted regular tree that is embedded in the graph.
 
To estimate the radius of the tree rooted at a given vertex, we must consider the asymptotic behavior of a function which is in fact Gauss' hypergeometric function $_{2}F_1\left(\begin{array}{c|}
1, \; 1-n \\
\frac{1-dn}{2} 
\end{array} ~   \frac{d}{2} \right)$.  While many asymptotic estimates for various hypergeometric functions are known, the particular one we need is not, and consequently we have the following corollary, 
\begin{co}\label{hyper}
$$_{2}F_1\left(\begin{array}{c|}
1, \; 1-n \\
\frac{1-dn}{2} 
\end{array} ~   \frac{d}{2} \right)
\propto \sqrt{n}$$
for all $d$ as $n \to \infty$.
\end{co}

\section{Results} 
First we note that for integer value random variable $E(X)=\sum nP(X=n)=\sum P(X \ge n)$ 
 
 We pick a vertex $v_0$ and now define a random variable $X$ to be the size i.e, the number of vertices in the embedded tree, of the largest $d$-regular tree rooted at the vertex $v_0$.  A tree is $d$ regular if all the vertices that are not leaves have a degree $d$.
The following lemma implies theorem \ref{main}

\begin{Lm}\label{technical}

For every $0\le \rho < \frac{1}{2}$ $$P(X \ge n^ \rho)\to 1$$

For  $\frac{1}{2}\le \rho < 1$ $$P(X \ge n^ \rho)\to 0$$ as $n\to \infty$

\end{Lm}
Note: $r \propto \log_d Ex(X)= \rho \log_d n$\\
We will use  a version of the configuration model. The configuration model, which is commonly use when working with random regular graphs, was developed mainly by Bolob\'as to construct a $d$-regular graph with $n$ vertices ($dn$ has to be even). Let $W=[n]\times [d]$.  A configuration is a partition  of $W$ into $nd/2$ pairs. These pairs are called edges. To construct a graph from a configuration, project $W$ into set of $V=[n]$ the vertices.  The vertices $i,j$ are joined by an edge iff $i$ and $j$ appear in an edge of the configuration $W$. The resulting graph might have loops and multiple edges, but since our motivation come from Belyi surfaces we will allow  loops and multiple edges.  It is standard technique when working with the configuration model to translate results from configurations to simple graphs without  loops and multiple edges.

\begin{figure}[htbp] 
   \centering
   \includegraphics[width=4in]{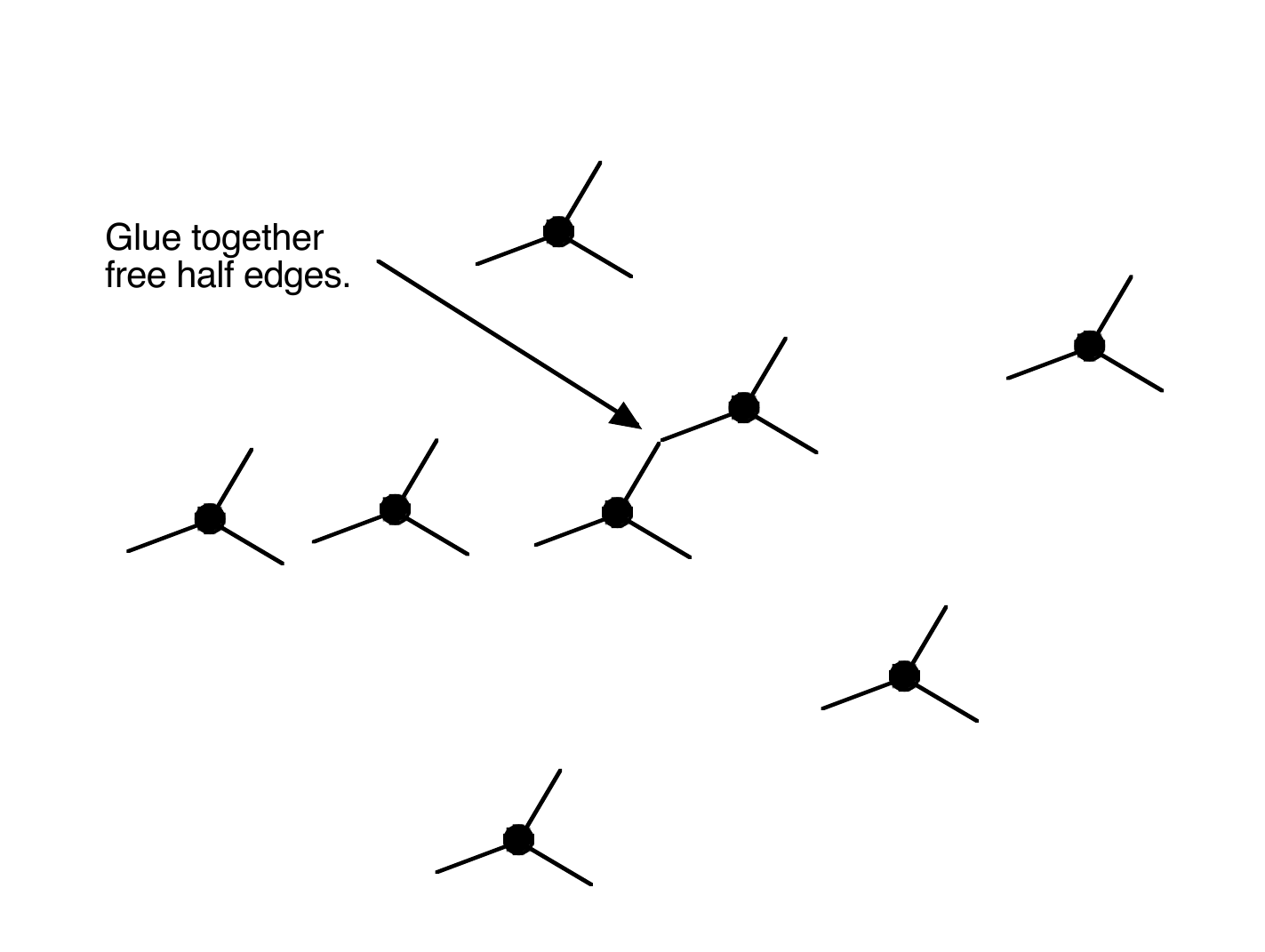} 
   \caption{The configuration model with $d=3$.}
   \label{config}
\end{figure}

We will use a variation of this model that was used be Brooks and Makover \cite{BM} for cubic graphs. 
We look at $W$ as a labeling of the set of $n$ vertices each of which has $d$ half edges attached to it (see Figure \ref{config}).   We will expose the pairs of the partition in "shells"  starting from the root, and at each stage we look at the free half edges that are attached  to the connected component around the root, identifying the second half of the edges connecting to the new shell of vertices around the root. We stop when both half edges belong to our connected component and therefore a loop has been created.

First we pick a vertex $v_0$, and note that $P(X\ge 1)=1$. 
Now we choose one of the half edges around $v_0$, and pick at random another half edge to glue to the our starting edge. There are $dn-1$ open half edges, for $d-1$ of them we will close a cycle, and therefore the processes of growing the tree will stop and $X=1$, For any other choice of half edge from the $dn-d$ "free" half edges we will get a tree of size $2$ and continue the process of growing the tree hence $X \ge 2$  Therefore $P(X\ge 2)=\frac{dn-d}{dn-1}$. We will continue by choosing the next half edge around $v_0$ from the remaining $dn-2$ open half-edges $2d-3$ will close a loop resulting in a tree of size $2$ and $nd-2d$ are "free" and the process continue therefore.
$$P(X\ge 3)=\frac{dn-d}{dn-1}\frac{dn-2d}{dn-3}$$
We observe that this pattern continues,
$$P(X\ge 4)=\frac{dn-d}{dn-1}\frac{dn-2d}{dn-3}\frac{dn-3d}{dn-5}$$
and so on.

We will start with the free half edges emanating from $v_0$, when we exhausted  all the half edges emanating from $v_0$, we will pick half edges emanating from the vertices that are of distance $1$ from $v_0$. Then distance $2$, $3$,  etc. In this way the tree is growing radially around $v_0$.
In general 

$$P(X\ge k)=\prod_{i=1}^{k-1}\frac{dn-id}{dn-(2i-1)}=\frac{d^{k-1}(n-1)!}{(n-k)!}\frac{ (dn-(2k-1))!!}{(dn-1)!!   }$$
where for the graphs under consideration, $n$ is even.
Next we will use Stirling's formula to get 

$$P(X\ge k) \propto \frac{d^{k-1}\sqrt{2\pi}(n-1)^{n-\frac{1}{2}} e^{-(n-1)}}{\sqrt{2\pi}(n-k)^{n-k+\frac{1}{2}} e^{-(n-k)}}
\frac{\sqrt{\pi} (dn-2k+1)^{\frac{dn-2k+2}{2}}e^{-\frac{dn-2k+1}{2}}}{{\sqrt{\pi} (dn-1)^{\frac{dn}{2}}e^{-\frac{dn-1}{2}}}}
$$

$$=\frac{(\frac{n-1}{n-k})^{n-\frac{1}{2}}}{(\frac{dn-1}{dn-2k+1})^\frac{dn}{2}}\left(\frac{dn-dk}{dn-2k+1}\right)^{k-1}$$
Now let $k=n^\rho$ for $0\le \rho \le 1$.  We need to consider the following limit:

$$\lim_{n\to \infty}\left(\frac{n-1}{n-n^\rho}\right)^{n-1}   \left (\frac{dn-2n^\rho+1}{dn-1}\right)^\frac{dn}{2}     \left(\frac{dn-dn^\rho}{dn-2n^\rho+1}\right)^{n^\rho-1}      $$
A tedious  but straightforward of l'H\^{o}pital's rule will show that

$$\lim_{n\to \infty}\left(\frac{n-1}{n-n^\rho}\right)^{n-1}   \left (\frac{dn-2n^\rho+1}{dn-1}\right)^\frac{dn}{2}     \left(\frac{dn-dn^\rho}{dn-2n^\rho+1}\right)^{n^\rho-1}  =  \left\{ \begin{array}{ll} 1 & \mbox{ } \rho< \frac{1}{2} \\ 

e^{-\frac{d-2}{2d}} &\mbox{   } \rho=\frac{1}{2} \\
0 &\mbox{ }\rho > \frac{1}{2} \end{array} \right. 
 \blacksquare $$
 
 To conclude the the proof of theorem \ref{main}  we need to estimate the rate of convergence of the tail of the distribution. 
 
 From \ref{technical}  we see that $$P(X \ge x n^{1/2})) \to \exp(-x^2\frac{d-2}{2d})$$ as $n\to \infty$, if $X$ is scaled by $n^{1/2}$ then $X$ as a distribution converges to a random variable with tail $$1-F(x)=\exp(-x^2\frac{d-2}{2d}).$$   This completes the proof of the theorem.

Notice that 
$$_{2}F_1\left(\begin{array}{c|}
1, \; 1-n \\
\frac{1-dn}{2} 
\end{array} ~   \frac{d}{2} \right)
=\sum_{k=1}^n\prod_{i=1}^{k-1}\frac{dn-id}{dn-(2i-1)}=\sum_{k=1}^n P(X\ge k)=E(X)
$$
where $_{2}F_1\left(\begin{array}{c|}
a, \; b \\
c
\end{array} ~  z \right)$ is Gauss's Hypergeometric function. 
This gives  an interesting corollary of Theorem \ref{main} about the asymptotic behavior of  Gauss's Hypergeometric function.

\begin{co}\label{hyper}
$$_{2}F_1\left(\begin{array}{c|}
1, \; 1-n \\
\frac{1-dn}{2} 
\end{array} ~   \frac{d}{2} \right)
\propto \sqrt{n}$$
for all $d$ as $n \to \infty$.
\end{co}

\bibliographystyle{amsplain}

 \end{document}